# Semilattices of Rectangular Bands and Groups of Order Two.

R. A. R. Monzo

**Abstract.** We prove that a semigroup *S* is a semilattice of rectangular bands and groups of order two if and only if it satisfies the identity $x = x^3$ and $xyx \in \{xy^2x, y^2x^2y\}$ $(x, y \in S)$.

## 1. Introduction

Semigroup theory has a certain symmetric elegance that links it to group theory. For example, a semigroup *S* is a group if (and only if) for every $x \in S$, $Sx = S = xS$. In addition, a semigroup *S* is a union of groups if (and only if) for every $x \in S$, $x \in Sx^2 \cap x^2S$ [3]. The powerful result that a semigroup is a union of groups if and only if it is a semilattice of completely simple semigroups is so well known that its beauty can almost be overlooked [3]. In this paper we apply this result to prove that the collection of all semigroups that are semilattices of semigroups that are either rectangular bands or groups of order two is a semigroup inclusion class [4]. Precisely, a semigroup is a semilattice of rectangular bands and groups of order two if and only if it satisfies the identity $x = x^3$ and $xyx \in \{xy^2x, y^2x^2y\}$ $(x, y \in S)$.

## 2. Notation, definitions and preliminary results

**Definition.** Let *G* be a group and $\Lambda \times I$ a non-empty set. Let $P: \Lambda \times I \to G$, with $P(\mu, j)$ denoted by $p_{\mu j}$. Let $S = I \times G \times \Lambda$ and define a product on S as follows: $(i, a, \mu)(j, b, \lambda) = (i, ap_{\mu j}b, \lambda)$. Then S endowed with this product is called a ***Rees*** $I \times \Lambda$ ***matrix semigroup (over the group* G *with the sandwich matrix P).***

The well known definition above is repeated here because it will be used extensively throughout this paper. All other terminology and notation can be found in [3]. Other well known results that will be used here follow.

**Result 1**. [5**,** Corollary IV.2.8] *A semigroup is completely simple if and only if it is isomorphic to a Rees matrix semigroup.*

**Result 2**. [3,Theorem 4.3] *A semigroup* S *is a union of groups if and only if for every* $x \in S$, $x \in x^2S \cap Sx^2$.

**Result 3.** [3, Theorem 4.6] *A semigroup is a union of groups if and only if it is a semilattice Y of completely simple semigroups* $S_\alpha = I_\alpha \times G_\alpha \times \Lambda_\alpha$ $(\alpha \in Y)$.

**Result 4.** [2, Proposition 1] *A semilattice Y of completely simple semigroups is a semilattice of rectangular groups if and only if the product of two idempotents is an idempotent.*

**Definition.** We will denote the collection of all left zero, right zero and rectangular bands by ***L*** $_0$, ***R*** $_0$ and ***RB*** respectively and that of all groups of order two by ***G*** $_2$. If S is a union of groups and $x \in S$ then $\mathbf{1}_x$ will denote the identity element of any group to which *x* belongs. If G is a group then the identity element of G is denoted by **1** or $\mathbf{1}_G$. Also if $x$ is an element of a semigroup S then $x^{-1}$ or $\left(x^{-1}\right)_G$ denotes the inverse of $x$ in any subgroup G of S. A ***rectangular group*** is the direct product of a rectangular band and a group.

Note that if an element $x$ of a semigroup belongs to two groups, G and H, then
$1_G = \left(x^{-1}\right)_G x = \left(x^{-1}\right)_G x 1_H = 1_G 1_H = 1_G x\left(x^{-1}\right)_H = x\left(x^{-1}\right)_H = 1_H$, and therefore $1_x$ is well-defined.
Similarly, $\left(x^{-1}\right)_G = 1_G \left(x^{-1}\right)_G = 1_H \left(x^{-1}\right)_G = \left(x^{-1}\right)_H x\left(x^{-1}\right)_G = \left(x^{-1}\right)_H 1_G = \left(x^{-1}\right)_H 1_H = \left(x^{-1}\right)_H$ and so $x^{-1}$ is well-defined.

**Definition.** [4] An ***inclusion class of semigroups*** is a collection of all semigroups that satisfy a given set of $k$ number of inclusions as follows:

(1) $W_1 = \{w_{1,1}, w_{1,2}, ..., w_{1,n_1}\} \subseteq \{t_{1,1}, t_{1,2}, ..., t_{1,m_1}\} = T_1$;

(2) $W_2 = \{w_{2,1}, w_{2,2}, ..., w_{2,n_2}\} \subseteq \{t_{2,1}, t_{2,2}, ..., t_{2,m_2}\} = T_2$; …;

(k) $W_k = \{w_{k,1}, w_{k,2}, ..., w_{k,n_k}\} \subseteq \{t_{k,1}, t_{k,2}, ..., t_{k,m_k}\} = T_k$, where the $w$'s and the $t$'s are semigroup words over some alphabet.

**Notation.** We write $[W_1 \subseteq T_1;\ W_2 \subseteq T_2;...;W_k \subseteq T_k]$ to denote the inclusion class determined by the set $\{W_i \subseteq T_i\}(i \in \{1,2,...,k\})$ of inclusions. If S is a semilattice and $\{\alpha, \beta\} \subseteq S$ then we write $\alpha \leq \beta$ if $\alpha = \alpha\beta$ and $\alpha < \beta$ if $\alpha = \alpha\beta$ and $\alpha \neq \beta$.

**Definition.** A semigroup S is a ***semilattice*** if $S \in [\ x = x^2;\ xy = yx\ ]$. A semilattice S is a ***chain*** if $\{\alpha, \beta\} \subseteq S$ and $\alpha \neq \beta$ implies $\alpha\beta \in \{\alpha, \beta\}$.

**Definition.** A semigroup S is ***a semilattice Y of semigroups*** $S_\alpha (\alpha \in Y)$ if S is a disjoint union of the $S_\alpha (\alpha \in Y)$ and for every $\{\alpha, \beta\} \subseteq Y$, $S_\alpha S_\beta \subseteq S_{\alpha\beta}$.

## 3. Some inclusion classes of semilattices of rectangular bands and groups of order two.

**Theorem 1.** *The following statements are equivalent:*

1. $S \in [\ xyx \in \{x, y\}\ ]$ *and*

2. S *is a chain Y of semigroups* $S_\alpha (\alpha \in Y)$ *where*

   *(1)* $S_\alpha \in RB \cup G_2 (\alpha \in Y)$,

   *(2) for any* $\alpha < \beta\ (\alpha, \beta \in Y)$ *and any* $x \in S_\alpha$, $y \in S_\beta$, $xy = yx = x$ *and*

   *(3) for any* $\alpha < \beta\ (\alpha, \beta \in Y)$ *with* $S_\alpha \in G_2$, $S_\alpha = \{1\}$.

Proof: $(1 \Rightarrow 2)$ Let $S \in [\, xyx \in \{x, y\}\,]$. Then clearly, $S \in [\, x = x^3 \,]$ and so by Result 2, S is a union of groups and, therefore, by Result 3, S is a semilattice $Y$ of completely simple semigroups $S_\alpha$ $(\alpha \in Y)$.

Let $x \in S_\alpha$ and $y \in S_\beta$ $(\alpha, \beta \in Y)$. Then either (a) $xyx = x$ and $yxy = y$ or (b) $xyx = x$ and $yxy = x$ or (c) $xyx = y$ and $yxy = y$ or (d) $xyx = y$ and $yxy = x$. In cases (a) and (d), $\alpha = \alpha\beta = \beta$. So if $\alpha \neq \beta$ then either case (b) or (c) holds. Therefore, $\alpha \neq \beta$ implies that either $\alpha\beta = \alpha$ or $\alpha\beta = \beta$. Hence, $Y$ is a chain.

We now show that $S_\alpha \in RB \cup G_2$ $(\alpha \in Y)$. Each $S_\alpha = I_\alpha \times G_\alpha \times \Lambda_\alpha$ $(\alpha \in Y)$. Assume that either $|I_\alpha| \neq 1$ or $|\Lambda_\alpha| \neq 1$. Then let $x = (i, g, \lambda)$ and $y = (j, h, \mu)$, where either $i \neq j$ or $\lambda \neq \mu$ and where $g$ and $h$ are arbitrary elements of $G_\alpha$. Now (b), (c) and (d) all imply that $i = j$ and $\lambda = \mu$. Thus, (a) holds.

Now $x = x^3$ and so $x^2 = 1_x$. So $(i, g, \lambda)^2 = (i, gp_{\lambda i}g, \lambda) = (i, (p_{\lambda i})^{-1}, \lambda)$. Hence, $gp_{\lambda i}g = (p_{\lambda i})^{-1}$. Let $g'$ be an arbitrary element of $G_\alpha$. Since $g$ was arbitrary we can let $g = g'(p_{\lambda i})^{-1}$. Then $(p_{\lambda i})^{-1} = g'(p_{\lambda i})^{-1} p_{\lambda i} g'(p_{\lambda i})^{-1}$ and hence $(g')^2$ is the identity element of $G_\alpha$ and therefore $G_\alpha \in G_2$. Hence, $G_\alpha$ is abelian and $p_{\lambda i} = (p_{\lambda i})^{-1}$.

Now (a) holds and so $xyx = x$. So $(i, g, \lambda) = (i, gp_{\lambda i}hp_{\mu i}g, \lambda)$ and therefore $g = gp_{\lambda i}hp_{\mu i}g$ and $1_x = p_{\lambda j}hp_{\mu i}g$. Setting $g = 1_{G_\alpha} = h$ gives $1_{G_\alpha} = p_{\lambda i}p_{\mu i}$. But $G_\alpha$ is abelian and so $g = gp_{\lambda i}hp_{\mu i}g = g^2 p_{\lambda i}p_{\mu i}h = h$. Therefore $G_\alpha = \{1\}$. This implies $S_\alpha \in RB$.

Now if $I_\alpha = \{i\}$ and $\Lambda_\alpha = \{\lambda\}$ then the mapping $g \mapsto (i, g(p_{\lambda i})^{-1}, \lambda)$ is an isomorphism between $G_\alpha$ and $S_\alpha$. As shown two paragraphs above, $G_\alpha \in G_2$ and so $S_\alpha \in G_2$. We have therefore shown that $S_\alpha \in RB \cup G_2$ $(\alpha \in Y)$ and so (1) is valid.

Now assume that $\alpha < \beta$ $(\alpha, \beta \in Y)$ and let $x \in S_\alpha$ and $y \in S_\beta$. Then (b) holds and so $xyx = x = yxy$. Then $1_x = x^{-1}x = x^{-1}yxy = 1_x 1_y$ and $1_x = xx^{-1} = yxyx^{-1} = 1_y 1_x$.

If $S_\alpha \in G_2$ then $S_\alpha$ is commutative and so
$xy = (xyx)y = [(xy)x]y = [x(xy)]y = x^2 y^2 = 1_x 1_y = 1_x = (yx)^2 = y(xyx) = yx$. However, if $z \in S_\alpha$ then $(xy)z = 1_x z = z = x(yz) = x1_z = x$. Therefore $S_\alpha = \{x\}$ and so $1_x = x = xy = yx$. (We have shown that (3) is valid.)

If $S_\alpha \in RB$ then $1_x = x = x^2 = (xyx)(yxy) = (xy)^3 = xy$. Also, from (b), $yx = (yx)^3 = (yxy)(xyx) = x^2 = x = xy$. Hence, $x = xy = yx$. We have shown that (2) is valid and this completes the proof of $(1 \Rightarrow 2)$.

$(2 \Rightarrow 1)$ Assume that S is a chain $Y$ of semigroups $S_\alpha$ $(\alpha \in Y)$ where (1), (2) and (3) in Theorem 1 are valid. We wish to show that $xyx \in \{x, y\}$ for any $\{x, y\} \subseteq S$. Since $S_\alpha \in RB \cup G_2$ $(\alpha \in Y)$, we can assume that $\alpha \neq \beta$.

Case 1. Suppose that $\alpha < \beta$ and $S_\alpha \in G_2$. Then by (3), $S_\alpha = \{x\}$ and so $xyx \in S_\alpha = \{x\}$.

Case 2. Suppose that $\alpha < \beta$ and $S_\alpha \in RB$. Then $xyx = (xy)x = (xy)x^2 = x(yx)x = x^2 = x$.

Case 3. If $\beta < \alpha$ and $S_\beta \in G_2$ then, by (3), $S_\beta = \{y\}$. Therefore, $xyx \in S_\beta = \{y\}$.

Case 4. As in the proof of case 2, $\beta < \alpha$ and $S_\beta \in RB$ implies $yxy = y$. Then

$$xyx = x(yxy)x = x(yx)^2 = xy^2 = (xy)y = y^2 = y.$$

We have therefore shown that $xyx \in \{x, y\}$, completing the proof of ($2 \Rightarrow 1$). *This completes the proof of Theorem 1.*

Note that we have shown that any $S \in [\, xyx \in \{x, y\}\,]$ is a chain $Y$ of rectangular bands, except possibly if $Y$ has a maximal element $\alpha$ and $S_\alpha \in G_2$ with $|S_\alpha| > 1$. The question arises as to whether the collection of chains of semigroups that are either rectangular bands or groups of order 2 is an inclusion class. In Theorem 5 below we prove this question in the affirmative when the word "chain" is replaced by "semilattice".

**Theorem 2.** *The following statements are equivalent:*

1. $S \in [\, xyx \in \{y, yx\}\,]$ and

2. $S$ *is a semilattice $Y$ of semigroups $S_\alpha (\alpha \in Y)$ where*

(1) $(S_\alpha)^2 \in R_0 \cup G_2 (\alpha \in Y)$,

(2) $\alpha < \beta$ and $(S_\alpha)^2 \in G_2$ implies $(S_\alpha)^2 = \{1\}$,

(3) $S_\alpha - (S_\alpha)^2 \neq \emptyset$ implies $(S_\alpha)^2 = \{1\}$,

(4) *for each $\alpha \in Y$ and $x \in S_\alpha$ there is a mapping $\theta_x : S \to S$ such that for any $\beta \in Y$,*
$\theta_x /_{S_\beta} : S_\beta \to (S_{\alpha\beta})^2$ *is a homomorphism, satisfying*

(4.1) *if* $(S_\alpha)^2 \in G_2$ *then for any* $g \in (S_\alpha)^2$ *and* $\{x, y\} \subseteq S_\alpha$, $\theta_g = \theta_{1_g} = \theta_x \theta_y = \theta_x$ *on* $S_\alpha$,

(4.2) *if* $\{\alpha, \beta, \gamma\} \subseteq Y$ *and* $(S_{\alpha\beta\gamma})^2 \in R_0$ *then* $\theta_{(\theta_x y)(\theta_y x)} = \theta_y \theta_x$ *on* $S_\gamma$ *and*

(4.3) *for every* $x \in S_\alpha$ *and* $y \in S_\beta$, $xy = (\theta_x y)(\theta_y x)$.

Proof: $(1 \Rightarrow 2)$ Assume that $S \in [\, xyx \in \{y, yx\}\,]$. First we will prove that $S \in [\, xy = (xy)^3\,]$. Let $\{x, y\} \subseteq S$. Then either (a) $xyx = y$ and $yxy = x$ or (b) $xyx = y$ and $yxy = xy$ or (c) $xyx = yx$ and $yxy = x$ or (d) $xyx = yx$ and $yxy = xy$. Note that since $S \in [\, xyx \in \{y, yx\}\,]$, $(xy)^3 \in \{xy, (xy)^2\}$. So in each case we can assume that $(xy)^3 = (xy)^2$.

Case (a): $xy = (yxy)(xyx) = (yx)^3 = [(xyx)(yxy)]^3 = (xy)^9 = (xy)^6 = (xy)^4 = (xy)^3$.

Case (b): $xy = y(xy) = (xyx)(yxy) = (xy)^3$.

Case (c): $(xy)^2 = (xyx)y = yxy = x = x(yxy) = x^2$. But then $x = (xy)^2 = x^2$ implies $x = (xy)^2 x$ and so $xy = \left[(xy)^2 x\right]y = (xy)^3$.

Case (d): $xy = yxy = (yx)y = (xyx)y = (xy)^2 = (xy)^3$.

By Results 2 and 3, $S^2$ is a semilattice $Y$ of completely simple semigroups $(S^2)_\alpha$ $(\alpha \in Y)$. We now show that the product of two idempotents of $S$ is an idempotent. Let $\{e, f\} \subseteq E_S$. Then $efe \in \{f, fe\}$. If $efe = f$ then $ef = e(efe) = efe = f \in E_S$. If $efe = fe$ then $(ef)^2 = (efe)f = fef \in \{e, ef\}$. So we can assume that $(ef)^2 = e$. Hence $ef = (ef)^2 f = (ef)^2$ and this completes the proof that the product of two idempotents is idempotent.

Now by Result 4, $S^2$ is a semilattice $Y$ of rectangular groups $(S^2)_\alpha$ $(\alpha \in Y)$. We want to show now that each $(S^2)_\alpha = L \times G \times R$ is either a right-zero semigroup or a group of order two.

Let $\{x, y\} \subseteq (S^2)_\alpha$ with $x = (i, g, \lambda)$ and $y = (j, h, \delta)$. Then $xyx \in \{y, yx\}$ and so $i = j$. So $(S^2)_\alpha = G \times R$. We have shown that an arbitrary rectangular group component of $S^2$ is a right group $G \times R$.

Note that, since we have already shown above that $xy = (xy)^3$, $G \in G_2$ and so $G$ is commutative and satisfies $g = g^{-1}$, $g^2 = 1$ and $ghg = h$ $(g, h \in G)$.

Now let $\{(g, r)(g, s)\} \subseteq G \times R$, with $r \neq s$. Then $(g, r) = (g, r)(g, s)(g, r) \in \{(g, s), (g^2, r)\}$. Therefore $g = g^2 = 1$. So either $|R| = 1$ or $|G| = 1$. Hence $G \times R$ is either a right-zero semigroup or a group of order 2.

For $\alpha \in Y$ we define $S_\alpha = \{x \in S : x^2 \in (S^2)_\alpha\}$. Let $\{x, y\} \subseteq S$ with $x^2 \in (S^2)_\alpha$ and $y^2 \in (S^2)_\beta$. We show that $xy \in (S^2)_{\alpha\beta} \subseteq S_{\alpha\beta}$. First note that $S$ is a null extension of a union of groups. From the proof of Theorem 5 [1], $xy \in xy^2 S \cap S x^2 y$. It is then straightforward to show that

(1) $xy = (xy1_{y^2})(xy1_{y^2})^{-1} xy(1_{x^2} xy)(1_{x^2} xy)^{-1}$ [*], with $1_{x^2} \in (S^2)_\alpha$ and $1_{y^2} \in (S^2)_\beta$.

Assume that $xy \in (S^2)_\delta$. Then, since $S^2$ is a semilattice of the semigroups $(S^2)_\alpha$ $(\alpha \in Y)$, it follows from (1) that $\delta = \delta\alpha\beta$.

Note that $x^3 \in (S^2)_\alpha$ and $y^3 \in (S^2)_\beta$. Therefore $x^3 y^3 = x^2(xy)y^2$, which implies that $\alpha\beta = \alpha\delta\beta = \delta$. Therefore $S$ is a semilattice of the semigroups $S_\alpha$ $(\alpha \in Y)$. It is straightforward to show that $(S_\alpha)^2 = (S^2)_\alpha$. Thus, we have proved part *(1)* of Theorem 2.

Suppose now that $\alpha < \beta$, $x \in (S^2)_\alpha \in G_2$ and $y \in (S^2)_\beta$. Then $xyx \in \{y, yx\}$ and so $xyx = yx$. But $yx \in (S^2)_{\alpha\beta} = (S^2)_\alpha \in G_2$ and so $x = 1$. Therefore $(S^2)_\alpha = \{1\}$. This proves part *(2)* of Theorem 2.

Suppose that $y \in S_\alpha - (S_\alpha)^2$. Let $x \in (S_\alpha)^2$. Then by hypothesis $xyx = yx$, and so $x = 1$. This proves part *(3)* of Theorem 2.

We proceed with the proof of part *(4)*. For any $x \in S_\alpha$ we define $\theta_x : S \to S$ as $y \mapsto xyx$ ($y \in S_\beta$). Note that $\theta_x /_{S_\beta} : S_\beta \to (S_{\alpha\beta})^2$, because $xyx = x(yx)^3 = [x(yx)^2](yx) \in (S^2)_{\alpha\beta} = (S_{\alpha\beta})^2$. Let $\{y, z\} \subseteq S_\beta$. If $(S_{\alpha\beta})^2 = (S^2)_{\alpha\beta} \in R_0$ then

$\theta_x(yz) = (xy)(zx) = xy(zx)^3 = (xyzx)(zx)^2 = (zx)^2 = [xyxx(zx)](zx)^2 = (xyx)[x(zx)^3] = (xyx)(xzx) = (\theta_x y)(\theta_x z)$

Suppose that $(S_{\alpha\beta})^2 \in G_2$. If $\alpha \neq \beta$ then either $\alpha\beta < \beta$ or $\alpha\beta < \beta$ and so, by *(2)*, $(S_{\alpha\beta})^2 = \{1\}$. This implies that $\{\theta_x(yz), \theta_x y, \theta_x z\} \subseteq (S_{\alpha\beta})^2 = \{1\}$ and so $\theta_x /_{S_\beta}$ is a homomorphism.

We can therefore assume that $\alpha = \beta = \alpha\beta$ and so $(S_{\alpha\beta})^2 = (S_\alpha)^2 = (S_\beta)^2 \in G_2$. Since, by hypothesis, $x^3 \in \{x, x^2\}$, $x^2 = x^4 = 1$. Then, $(\theta_x y)(\theta_x z) = xyx^2 zx = xyzx = \theta_x(yz)$ and so $\theta_x /_{S_\beta}$ is a homomorphism in any case. This proves *(4)*.

Let $(S_\alpha)^2 \in G_2$, $g \in (S_\alpha)^2$, $1 = 1_g$ and $\{x, y, z\} \subseteq S_\alpha$. Then $\theta_g x = gxg = g1x1g = g^2 1x1 = 1x1 = \theta_1 x$ and $\theta_x \theta_y z = xyzyx = yx(xyz) = yx^2 yz = yx^4 yz = y^2 z = y^4 z = 1z = 1z1 = \theta_1 z$. Also, $\theta_x z = xzx = x1zx = x1z1x = (1x)x1z = 1x^2 1z = 1z = x^4 z = x(x^2)(xz) = x(xz)x^2 = x^2 zx^2 = 1z1 = \theta_1 z$. Hence, $\theta_g = \theta_{1_g} = \theta_x \theta_y = \theta_x$ and this proves *(4.1)*.

Let $\{\alpha, \beta, \gamma\} \subseteq Y$ and $(S_{\alpha\beta\gamma})^2 \in R_0$. Suppose $x \in S_\alpha$, $y \in S_\beta$ and $z \in S_\gamma$. Then
$\theta_{(\theta_x y)(\theta_y x)} z = \theta_{(xyx)(yxy)} z = (xy)^3 z(xy)^3 = xyzxy = (xyzxy)^2 = (xyzxyxy)(zxy) = zxy = (yxzyx)(zxy) = $
$= (yxz)(yxzxy) = yxzxy = \theta_y \theta_x z$. This proves *(4.2)*.

Finally, *(4.3)* follows from the fact that $xy = (xy)^3$. This completes the proof of $(1 \Rightarrow 2)$.

$(2 \Rightarrow 1)$ Assume that the hypotheses of the "only if" part of Theorem 2 are valid. We first show that the product defined is associative. Let $x \in S_\alpha$, $y \in S_\beta$, $z \in S_\gamma$. Using *(4.3)* we need to show that:
$[\theta_{(\theta_x y)(\theta_y x)} z]\{\theta_z [(\theta_x y)(\theta_y x)]\} = \{\theta_x [(\theta_y z)(\theta_z y)]\}[\theta_{(\theta_y z)(\theta_z y)} x]$. However, since by *(4)* $\theta_z$ and $\theta_x$ are homomorphisms on $(S_{\alpha\beta})^2$ and $(S_{\beta\gamma})^2$ respectively, this equation becomes:

(2) $(\theta_{(\theta_x y)(\theta_y x)} z)(\theta_z \theta_x y)(\theta_z \theta_y x) = (\theta_x \theta_y z)(\theta_x \theta_z y)(\theta_{(\theta_y z)(\theta_z y)} x)$,

with – by *(4)* again -- each of the 6 terms an element of $(S_{\alpha\beta\gamma})^2$. If $(S_{\alpha\beta\gamma})^2 \in R_0$ then, since by hypothesis *(4.2 )*, $\theta_{(\theta_y z)(\theta_z y)} x = \theta_z \theta_y x$, equation (2) is valid in this case.

So we can assume that $(S_{\alpha\beta\gamma})^2 \in G_2$. We can therefore assume that $\alpha = \beta = \gamma$, or else there exists $\sigma \in \{\alpha, \beta, \gamma\}$ such that $\alpha\beta\gamma < \sigma$, which implies $(S_{\alpha\beta\gamma})^2 = \{1\}$. [This would imply that equation $(2)$ is valid.]. But if $\alpha = \beta = \gamma$ then, by *(4.1 )*, each side of equation $(2)$ equals, $(\theta_1 z)(\theta_1 y)(\theta_1 x)$, so $(2)$ is valid.

Now we need to prove that for $x \in S_\alpha$ and $y \in S_\beta$, $xyx \in \{y, yx\}$. If $(S_{\alpha\beta})^2 \in R_0$ then, since $\{xyx, xyxx\} \subseteq (S_{\alpha\beta})^2$, $xyx = (xyx)^2 = (xyxx)(yx) = yx$. We can assume, therefore, that $(S_{\alpha\beta})^2 \in G_2$.

If $\alpha \neq \beta$ then either $\alpha\beta < \alpha$ or $\alpha\beta < \beta$. In either case $\{xyx, yx\} \subseteq (S_{\alpha\beta})^2 = \{1\}$ and so $xyx = yx$. We can assume therefore that $\alpha = \beta$. We can also assume that $y \in (S_\alpha)^2$, or else by *(3)*, $\{xyx, yx\} \subseteq (S_\alpha)^2 = \{1\}$. Note that since $xy = (\theta_x y)(\theta_y x) \in (S_\alpha)^2 \in G_2$, $(S_\alpha)^2$ is an abelian group. Therefore, $xy = (\theta_x y)(\theta_y x) = (\theta_y x)(\theta_x y) = yx$. Also, $x^2 = (\theta_x x)^2 = 1$. Hence, $xyx = yx^2 = y1 = y$. This completes the proof that $S \in [\ xyx \in \{y, yx\}\ ]$.

*So the proof of Theorem 2 is complete.*

**Corollary 3**. *The following statements are equivalent:*

*1. S is a semilattice Y of semigroups $S_\alpha (\alpha \in Y)$ where*

    *(1.1) $S_\alpha \in R_0 \cup G_2 (\alpha \in Y)$,*

    *(1.2) $\{\alpha, \beta\} \subseteq Y$, $\alpha < \beta$ and $S_\alpha \in G_2$ implies $S_\alpha = \{1\}$*

    *(1.3) there is a collection of mappings $\{\theta_x : S \to S\ /\ x \in S\}$ satisfying the following properties:*

        *(a) for $\{\alpha, \beta\} \subseteq Y$ each $\theta_x /_{S_\beta} : S_\beta \to S_{\alpha\beta}$ is a homomorphism,*

        *(b) $\{\alpha, \beta, \gamma\} \subseteq Y$, $x \in S_\alpha$, $y \in S_\beta$ and $S_{\alpha\beta\gamma} \in R_0$ imply $\theta_{(\theta_x y)(\theta_y x)} = \theta_y \theta_x$ on $S_\gamma$ and*

        *(c) if $S_\alpha \in G_2$ and $\{x, y\} \subseteq S_\alpha$ then $\theta_1 = \theta_x \theta_y = \theta_x$ on $S_\alpha$ and*

*2. $S \in [\ xyx \in \{y, yx\};\ x = x^3\ ]$.*

**Theorem 4.** *The following statements are equivalent:*

1. S *is a semilattice Y of semigroups* $S_\alpha$ $(\alpha \in Y)$ *with* $S_\alpha \in R_0 \cup G_2$ *and*

2. $S \in [\, xyx \in \{yx, y^2x^2y\};\ x = x^3\,]$.

Proof: $(1 \Rightarrow 2)$ For any $x \in R_0 \cup G_2$, $x = x^3$ and therefore $S \in [\, x = x^3\,]$. Suppose that $x \in S_\alpha$ and $y \in S_\beta$ $(\alpha, \beta \in Y)$. If $S_{\alpha\beta} \in R_0$ then $xyx = (xyx)^2 = (xyx^2)yx = yx$. Suppose that $S_{\alpha\beta} \in G_2$. Then,

$xy = (xy)^3 = x(yxy)xy = x(xy)yxy = x^2y^2xy = x^2(y^21)xy = x^2xyy^21 = x^3y^3 = x^3yy^21 = y^21x^3y =$
$= y^2x^3y1 = y1y^2x^3 = y^3x^3 = yx$. Also, $1 = (1y)^2 = 1y1y = 1yy = 1y^2 = y^21$. Then,

$xyx = xxy = 1x^2y = (y^21)x^2y = y^2x^2y$. We have proved therefore that $S \in [\, xyx \in \{yx, y^2x^2y\}\,]$. Hence, $S \in [\, x = x^3\,] \cap [\, xyx \in \{yx, y^2x^2y\}\,]$.

$(2 \Rightarrow 1)$ If $\{e, f\} \subseteq E_S$ then, since $efe \in \{fe, fef\}$, $ef = (ef)^3 = ef(efe)f \in \{ef(fe)f, ef(fef)f\} = \{(ef)^2\}$. It then follows from Results 1, 2 and 3 that $S$ is a semilattice $Y$ of rectangular groups $S_\alpha$ $(\alpha \in Y)$. Let $S_\alpha = L_\alpha \times G_\alpha \times R_\alpha$. The fact that $S \in [\, xyx \in \{yx, y^2x^2y\}\,]$ implies $|L_\alpha| = 1$ and so $S_\alpha = G_\alpha \times R_\alpha$. Since $S \in [\, x = x^3\,]$, $G_\alpha \in [\, x = x^3\,]$ and $G_\alpha \in G_2$. Let $\{x = (g, r),\ y = (h, s)\} \subseteq S_\alpha$. Then $xyx = (ghg, r) = (h, r) \in \{(hg, r), (h^2g^2h, s)\} = \{(hg, r), (h, s)\}$. So $|R_\alpha| > 1$ implies $|G_\alpha| = 1$. Thus, either $|R_\alpha| = 1$ or $|G_\alpha| = 1$. So $S_\alpha \in R_0 \cup G_2$, which is what we needed to prove. *This completes the proof of Theorem 4.*

**Theorem 5**: *The following statements are equivalent*:

1. S *is a semilattice Y of semigroups* $S_\alpha$ $(\alpha \in Y)$ *with* $S_\alpha \in RB \cup G_2$ $(\alpha \in Y)$ *and*

2. $S \in [\, xyx \in \{xy^2x, y^2x^2y\}\ ;\ x = x^3\,]$.

Proof: $(1 \Rightarrow 2)$ Clearly, $x = x^3$ for any $x \in S$. Let $\{\alpha, \beta\} \subseteq Y$ with $x \in S_\alpha$ and $y \in S_\beta$. If $S_{\alpha\beta} \in RB$ then $xyx = (xyx)^3 = xy(x^2yx^2)yx = xyyx = xy^2x$.
So we can assume that $S_{\alpha\beta} \in G_2$. Then $(xy)^2 = 1 = xyx1y = x(1y)yx = x1y^2x = y^2x(x1) = y^2x^21 =$
$y^21(x^21) = x^21(y^21) = x^2(1y^21) = x^2(y^21) = (x^2y^2)1 = x^2y^2$. Therefore, $1 = x^2y^2 = y^2x^2$.
Then, $xyx = (xyx)^3 = x(yx^2)(yx^2y)x = x(yx^2y)(yx^2)x = x(yx)(xy^2)x^3 = x(xy^2)(yx)x^3 = x^2y^2(yx^4) =$
$= 1yx^4 = yx^4 = y^3x^4 = y(y^2x^2)x^21 = y1x^21 = yx^21 = yx^2 = y^3x^2 = y(y^2x^2) = y1 = 1y1 = 1y = (y^2x^2)y$.
Hence, $S \in [\, xyx \in \{xy^2x, y^2x^2y\}\ ;\ x = x^3\,]$.

$(2 \Rightarrow 1)$ Since $S \in [\, x = x^3\,]$, by Results 2 and 3, S is a semilattice of completely simple semigroups.

We now show that the product $ef$ of the idempotents $e$ and $f$ is idempotent. We have

$efefe = efe(fef)(efe) \in \{efe(fef)^2 efe, (fef)^2(efe)^2 fef\} = \{(ef)^5 e, (fe)^6 f\} = \{efe, (fe)^2 f\}$.

If $efefe = efe$ then $(ef)^2 = (efe)f = (efefe)f = (ef)^3 = ef$. So we can assume that $efefe = (fe)^2 f$.

Then, $efefe = (fe)^2 f = [(fe)^2 f]e = (fe)^3 = fe = [(fe)^2 f]f = (efefe)f = (ef)^3 = ef$.

Therefore, $(ef)^2 = e(fe)f = e(ef)f = ef$. So by Result 3, S is a semilattice $Y$ of rectangular groups $S_\alpha = G_\alpha \times E_\alpha (\alpha \in Y)$, where $E_\alpha \in RB(\alpha \in Y)$.

Fix $\alpha \in Y$ and let $\{x,y\} \subseteq S_\alpha$ with $x = (1,(i,\lambda))$ and $y = (h,(j,\sigma))$ where $h$ is an arbitrary element of $G_\alpha$ and $(i,\lambda)$ and $(j,\sigma)$ are arbitrary elements of $E_\alpha$. Then,

$xyx = (h,(i,\lambda)) \in \{xy^2 x, y^2 x^2 y\} = \{(h^2,(i,\lambda)),(h,(j,\sigma))\}$. So either $E_\alpha$ has only one element or $h = h^2$ for every $h \in G_\alpha$. Hence $S_\alpha$ is isomorphic to $G_\alpha$ or to $E_\alpha$. In the former case, since $S \in [x = x^3]$, $G_\alpha \in G_2$.
*This completes the proof of Theorem 5.*
**References**

[1] Clarke, G.T., Monzo, R.A.R.,: *A Generalisation of the Concept of an Inflation of a Semigroup.* Semigroup Forum **60,** 172-186 (2000)

[2] Clifford, A.H.: *The Structure of Orthodox Unions of Groups.* Semigroup Forum **3,** 283-337 (1971).

[3] Clifford, A.H., Preston, G.B.: The Algebraic Theory of Semigroups. Math. Surveys of the American Math Soc., vol. 1. Am. Math. Soc., Providence (1961).

[4] Monzo, R.A.R.: *Further results in the theory of generalised inflations of semigroups.* Semigroup Forum **76** 540-560 (2008).

[5] Petrich, M.: Introduction to Semigroups. Charles E. Merrill Publishing Company. Columbus, Ohio. (1973).